\documentclass[12pt]{amsart}

\usepackage{amsthm}

\usepackage{mathrsfs,mathtools,latexsym,eucal}
\usepackage{amscd,amsfonts,amssymb,amsmath,amsthm}
\usepackage{mathrsfs}
\usepackage{graphicx,graphics,color}
\usepackage[utf8]{inputenc}
\usepackage[active]{srcltx}
\usepackage[english]{babel}
\usepackage[pagewise,displaymath,mathlines]{lineno}
\usepackage{listings}
\usepackage{color}
\usepackage{subfigure} 
\usepackage{url}
\usepackage{float}
\usepackage{enumerate}

\usepackage[active]{srcltx}

\usepackage{pgf,tikz}
\usetikzlibrary{arrows}

\definecolor{dkgreen}{rgb}{0,0.6,0}
\definecolor{gray}{rgb}{0.5,0.5,0.5}
\definecolor{mauve}{rgb}{0.58,0,0.82}

\newtheorem*{theorem*}{Theorem}
\newtheorem{theorem}{Theorem}[section]

\newtheorem{corollary}[theorem]{Corollary}

\newtheorem{definition}[theorem]{Definition}

\newcommand{\R}{\mathbb{R}}
\newcommand{\N}{\mathbb{N}}

\renewcommand{\leq}{\leqslant}
\renewcommand{\geq}{\geqslant}

\numberwithin{equation}{section}

\newcommand{\card}{\mathrm{card}}

\newcommand{\disp}{\displaystyle}
\newcommand{\vareps}{\varepsilon}

\begin{document}

\title{Characterizations of a Banach space through the strong lacunary and the lacunary statistical summabilities}

\author{Soledad Moreno-Pulido}
\address{Department of Mathematics, College of Engineering, University of Cadiz, Puerto Real 11510, Spain (EU)}
\email{{\tt soledad.moreno@uca.es}}

\author{Giuseppina Barbieri}
\address{Department of Mathematics,
                   University of Salerno,
                   Via Giovanni Paolo II,
                   84084 Fisciano,
                   ITALY}
\email{{\tt gibarbieri@unisa.it}}

\author{Fernando León-Saavedra}
\address{Departamento de Matemáticas, Universidad de Cádiz, Facultad de Ciencias Sociales y de la Comunicación, Avenida de la Universidad s/n, 11403-Jerez de la Frontera (Cádiz)}
\email{{\tt fernando.leon@uca.es}}

\author{Francisco Javier Pérez-Fernández}
\address{Departamento de Matemáticas, Universidad de Cádiz, Facultad de Ciencias, Avenida de la Universidad s/n, Puerto Real (Cádiz)}
\email{{\tt javier.perez@uca.es}}

\author{Antonio Sala-Pérez}
\address{Department of Mathematics, College of Engineering, University of Cadiz, Puerto Real 11510, Spain (EU)}
\email{{\tt antonio.sala@uca.es}}

 \begin{abstract}
 In this paper we characterize the completeness of a normed space through the strong lacunary ($N_\theta$) and lacunary statistical convergence ($S_\theta$) of series. A new characterization of weakly unconditionally Cauchy series and unconditionally convergent series through $N_\theta$ and $S_\theta$ is obtained. We also relate the summability spaces associated with these summabilities with the strong $p$-Ces\`aro convergence summability space.
 \end{abstract}


 \maketitle \thispagestyle{empty}

 \section{Introduction}

Let $X$ be a normed space, a sequence $(x_k)\subset X$ is said to be \emph{strongly 1-Ces\`{a}ro summable  (briefly, $|\sigma_1|$-summable) to $L\in X$} if
$$\lim_{n\to \infty} \frac{1}{n}\sum_{k=1}^n \|x_k-L\|=0.$$
This type of summability was introduced by Hardy-Littlewood \cite{hardy} and Fekete \cite{fekete} and it is related to the convergence of Fourier series (see \cite{boos,zeller-beekmann}). The $|\sigma_1|$ summability along with the statistical convergence \cite{zygmund} started a very striking theory with important applications \cite{JIneq3,JIneq5,MursaleenLibro}. Some years later, the \emph{strong lacunary summability $N_\theta$} was presented by Freedman et al. \cite{Freedman78} by introducing lacunary sequences and showed that $N_\theta$ is a larger class of $BK$-spaces which had many of the characteristics of $|\sigma_1|$. Later on, Fridy \cite{Fridy02,Fridy01} showed the concept of statistical lacunary summability and they related it with the statistical convergence and the $N_\theta$ summability.

The characterization of a Banach space through different types of convergence  has been dealt by authors like Kolk \cite{kolk}, Connor, Ganichev and Kadets \cite{connorganichevkadets},\dots

Let  $\sum x_i$ be  a series in  a normed space $X$,  in \cite{Aizpuru2000} the authors introduced the space of convergence $S(\sum x_i)$ associated to the series $\sum x_i$, it is defined as the space of sequences $(a_j)$ in $\ell_\infty$ such that $\sum a_ix_i$ converges. They also prove that the space $X$ is complete if and only if for every weakly unconditionally Cauchy series $\sum x_i$, the space $S(\sum x_i)$ is complete. Recall that a series is called weakly unconditionally Cauchy (wuC) if for every permutation $\pi$ of $\N$, the sequence $(\sum_{i=1}^n x_{\pi(i)})$ is a weakly Cauchy sequence. We will also rely in a powerful known result that states that a series $\sum x_i$ is wuC if and only if $\sum |f(x_i)|<\infty$ for all $f\in X^\ast$ (see \cite{Diestel1984} for Diestel's complete monograph about series in Banach spaces). 

In \cite{SalaFilomat,Leon-Saavedra2019} a Banach space is characterized by means of the strong $p$-Ces\`{a}ro summability ($\textrm{w}_p$) and ideal-convergence. In this manuscript, the $N_\theta$ and $S_\theta$ summabilities are used along with the concept of weakly unconditionally series to characterize a Banach space. In Section 2 we introduce these two kinds of summabilities which are regular methods and we recall some properties. In Section 3 and 4 we introduce the spaces $S_{S_{\theta}}(\sum_i x_i)$ and $S_{S_{\theta}}(\sum_i x_i)$ which will be used in Section 5 to characterize the completeness of a space.

\section{Preliminaries}

In this section we present the definition of $N_\theta$ and $S_\theta$ summabilities for Banach spaces and the relations between them. First, we recall the concept of lacunary sequences.

\begin{definition}
A \emph{lacunary sequence} is an increasing sequence of natural numbers $\theta=(k_r)$ such that $k_0=0$ and $h_r=k_r-k_{r-1}\to+\infty$ as $r\to\infty$. The intervals determined by $\theta$ will be denoted by $I_r=(k_{r-1},k_r]$ and the ratio $\displaystyle\frac{k_r}{k_{r-1}}$ will be denoted  by $q_r$.
\end{definition}

We now give the definition of strong lacunary summability for Banach spaces based on the one given by Freedman for real-valued sequences \cite{Freedman78}.

\begin{definition}
Let $X$ be a Banach space and $\theta=(k_r)$ a lacunary sequence. A sequence $x=(x_k)$ in $X$ is \emph{lacunary strongly convergent} or $N_{\theta}-$summable to $L\in X$ if $\displaystyle\lim_{r\to\infty}\frac{1}{h_r}\sum_{k\in I_r}\|x_k-L\|=0$,  and we write $N_{\theta}$-$\lim x_k=L$ or $x_k\underset{N_{\theta}}{\to} L$.
\end{definition}

Let $N_{\theta}$ be the space of all lacunary strongly convergent sequences,
$$N_{\theta}=\left\{(x_k)\subseteq X:\lim_{r\to\infty}\frac{1}{h_r}\sum_{k\in I_r}\|x_k-L\|=0\mbox{ for some }L\right\}.$$
The space $N_{\theta}$ is a BK$-$space endowed with the norm $\|x_k\|_{\theta}=\displaystyle\sup_r\displaystyle\frac{1}{h_r}\sum_{k\in I_r}\|x_k\|$.

In 1993, Fridy and Orhan \cite{Fridy01} introduced a generalization of the statistical convergence, the lacunary statistical convergence, using lacunary sequences. To accomplish this, they substituted the set $\{k:k\leq n\}$ by the set $\{k:k_{r-1}< k\leq k_r\}$. We recall now the definition of $\theta-$density of a subset $K\subset\N$.

\begin{definition}
Let $\theta=(k_r)$ be a lacunary sequence. If $K\subset\mathbb{N}$, the $\theta-$density of $K$ is denoted by
     $d_{\theta}(K)=\displaystyle\lim_r\frac{1}{h_r}\card(\{k\in I_r: k\in K\}),$
whenever this limit exists.
\end{definition}

It is easy to show that this density is a finitely additive measure and we can define the concept of lacunary statistically convergent sequences for Banach spaces. 

\begin{definition}
Let $X$ be a Banach space and $\theta=(k_r)$ a lacunary sequence. A sequence $x=(x_k)$ is a \emph{lacunary statistically convergent sequence} to $L\in X$ if given $\varepsilon>0$,
$$d_{\theta}(\{k\in I_r: \|x_k-L\|\geq\varepsilon\})=0,$$
or equivalently,
$$d_{\theta}(\{k\in I_r: \|x_k-L\|<\varepsilon\})=1,$$
we say that $(x_k)$ is $S_{\theta}$-convergent and we write $x_k\to_{S_{\theta}}L$.\end{definition}

\begin{theorem}
Let $X$ be a Banach space and $(x_k)$ a sequence in $X$. Notice that $S_\theta$ and $N_\theta$ are regular methods.

{\it Proof.}

\begin{enumerate}
    \item If $(x_k)\to L$, then $(x_k)\underset{N_\theta}{\to} L$.

        Let $\varepsilon>0$, then  there exists $k_0$ such that if $k\geq k_0$, then
        $$
            \|x_k-L\|<\varepsilon.
        $$
        Hence  there exists $r_0\in\mathbb{N}$ with  $r_0\geq k_0$ such that  if $r\geq r_0$ we have
        \begin{align*}
            \frac{1}{h_r}\sum_{k\in I_r}\|x_k-L\|<\frac{1}{h_r}\sum_{k\in I_r}\varepsilon=\frac{h_r}{h_r}\varepsilon=\varepsilon
        \end{align*}
        which implies that $\displaystyle\lim_{r\to\infty}\frac{1}{h_r}\sum_{k\in I_r}\|x_k-L\|=0$.
    \item If $(x_k)\to L$, then $(x_k)\underset{S_\theta}{\to} L$.

        Simply observe that, since $(x_k)\to L$, given $\varepsilon>0$ there exists $k_0$ such that for every $k\geq k_0$ we get $\card(\{k\in I_r: \|x_k-L\|\geq\varepsilon\})=0$, which implies $d_\theta(\{k\in I_r: \|x_k-L\|\geq\varepsilon\})=0$ for every $k\geq k_0$.
\end{enumerate}
\end{theorem}

Fridy and Orhan \cite{Fridy02} showed that $N_{\theta}$ and $S_{\theta}$ are equivalent for real-valued bounded sequences. This fact also holds for Banach spaces and we include the proof for the sake of completeness.

\begin{theorem}\label{FridyAcotadas}
    Let X be a Banach space, $(x_k)$ a sequence in $X$ and $\theta=(k_r)$ a lacunary sequence. Then:
    \begin{enumerate}
        \item $(x_k)\underset{N_\theta}{\to} L$ implies  $(x_k)\underset{S_\theta}{\to} L$.
        \item $(x_k)$ bounded and $(x_k)\underset{S_\theta}{\to} L$ imply  $(x_k)\underset{N_\theta}{\to} L$.
    \end{enumerate}
\end{theorem}

{\it Proof.} 1. If $(x_k)\underset{N_\theta}{\to} L$, then for every $\vareps>0$, 
$$
    \sum_{k\in I_r}\|x_k-L\|\geq \sum_{\substack{k\in I_r\\ \|x_k-L\|\geq\varepsilon}}\|x_k-L\|\geq\vareps\;\card(\{k\in I_r:\|x_k-L\|\geq\vareps\}),
$$
which implies that $(x_k)\underset{S_\theta}{\to} L$.

2. Let us suppose that $(x_k)$ is bounded and $(x_k)\underset{S_\theta}{\to} L$. Since $(x_k)$ is bounded, there exists $M>0$ such that $\|x_k-L\|\leq M$ for every $k\in\N$. Given $\vareps>0$,
\begin{align*}
    \frac{1}{h_r}\sum_{k\in I_r}\|x_k-L\| & =\frac{1}{h_r}\sum_{\substack{k\in I_r\\ \|x_k-L\|\geq\varepsilon}}\|x_k-L\|+\frac{1}{h_r}\sum_{\substack{k\in I_r\\ \|x_k-L\|<\varepsilon}}\|x_k-L\| \\
    & \leq \frac{M}{h_r}\card(\{k\in I_r: \|x_k-L\|\geq\vareps\})+\vareps,
\end{align*}
so we deduce that $(x_k)\underset{N_\theta}{\to} L$.

We now give the definition of lacunary statistically Cauchy sequences in Banach spaces as a generalization of the definition for real-valued sequences by Fridy and Orhan in \cite{Fridy01}.

\begin{definition}
Let $X$ be a Banach space and $\theta=(k_r)$ a lacunary sequence. A sequence $x=(x_k)$ is a \emph{lacunary statistically Cauchy} sequence if there exists a subsequence $x_{k'(r)}$ of $x_k$ such that $k'(r)\in I_r$ for every $r\in\N$, $\displaystyle\lim_{r\to \infty}x_{k'(r)}=L$ for some $L\in X$ and for every $\varepsilon>0$,
$$\lim_{r\to\infty}\frac{1}{h_r}\card(\{k\in I_r:\|x_k-x_{k'(r)}\|\geq\varepsilon\})=0,$$
or equivalently,
$$\lim_{r\to\infty}\frac{1}{h_r}\card(\{k\in I_r:\|x_k-x_{k'(r)}\|<\varepsilon\})=1.$$
In this case we say that $(x_k)$ is  $S_{\theta}$-Cauchy.\end{definition}

An important result in \cite{Fridy01} is the $S_{\theta}$-Cauchy Criterion and some of the next theorems in this work rely on it. This result can also be obtained for sequences in Banach spaces, and we include the proof for the sake of completeness.

\begin{theorem}\label{FridyCauchy}
    Let $X$ be a Banach space. A sequence $(x_k)$ in $X$ is $S_\theta$-convergent if and only if it is $S_{\theta}$-Cauchy.
\end{theorem}
    
{\it Proof.} Let $(x_k)$ be an $S_\theta$-convergent sequence in X and for every $k\in\N$, we define $K_j=\{k\in\N:\|x_k-L\|<1/j\}$. Observe that $K_j\supseteq K_{j+1}$ and $\disp\frac{\card(K_j\cap I_r)}{h_r}\to 1$ as $r\to\infty$.

Set $m_1$ such that if $r\leq m_1$ then $\card(K_1\cap I_r)/h_r>0$, that is, $K_1\cap I_r\neq\varnothing$. Next, choose $m_2>m_1$ such that if $r\geq m_2$, then $K_2\cap I_r\neq\varnothing$. Now, for each $m_1\leq r\leq m_2$, we choose $k'_r\in I_r$ such that $k'_r\in I_r\cap K_1$, i.e., $\|x_{k'_r}-L\|<1$. Inductively, we choose $m_{p+1}>m_p$ such that if $r>m_{p+1}$, then $I_r\cap K_{p+1}\neq\varnothing$. Thus, for all $r$ such that $m_p\leq r<m_{p+1}$, we choose $k'_r\in I_r\cap K_p$, and we have $\|x_{k'_r}-L\|<1/p$.

Therefore, we have a sequence $k'_r$ such that $k'_r\in I_r$ for every $r\in\N$ and $\lim_{r\to\infty}x_{k'_r}=L$. Finally,
\begin{align*}
    \frac{1}{h_r}\card(\{k\in I_r:\|x_k-x_{k'_r}\|\geq\vareps\})\leq& \frac{1}{h_r}\card(\{k\in I_r:\|x_k-L\|\geq \vareps/2\})\\
               &+\frac{1}{h_r}\card(\{k\in I_r:\|x_{k'_r}-L\|\geq \vareps/2\}).
\end{align*}
Since $(x_k)\underset{S_\theta}{\to} L$ and $\lim_{r\to\infty}x_{k'_r}=L$ we deduce that $(x_k)$ is $S_\theta$- Cauchy.

Conversely, if $(x_k)$ is a Cauchy sequence, for every $\vareps>0$,
\begin{align*}
    \card(\{k\in I_r:\|x_k-L\|\}\|\geq\vareps\})\leq& \card(\{k\in I_r:\|x_k-x_{k'_r}\|\geq \vareps/2\})\\
               &+\card(\{k\in I_r:\|x_{k'_r}-L\|\geq \vareps/2\}).
\end{align*}

Since $(x_k)$ is $S_\theta$-Cauchy and $\lim_{r\to\infty}x_{k'_r}=L$, we deduce that $(x_k)\underset{S_\theta}{\to} L$.

\section{The statistical lacunary summability space}

Let $\sum_i x_i$ be a series in a real Banach space $X$ and $\theta=(k_r)$ a lacunary sequence. We define
$$
	S_{S_{\theta}}\left(\sum_i x_i\right)=\left\{(a_i)_i\in\ell_{\infty}: \sum_{i} a_i x_i\mbox{ is }S_{\theta}\mbox{-summable}\right\}
$$
endowed with the supremum norm. This space will be called the space of $S_{\theta}$-summability associated to the series $\sum_i x_i$. The following theorem characterizes the completeness of the space $S_{S_{\theta}}\bigl(\sum_i x_i\bigr)$.

\begin{theorem}\label{S_tcompleto}
	Let $X$ be a real Banach space and $\theta=(k_r)$ a lacunary sequence. The following conditions are equivalent:
    \medskip
    \begin{enumerate}[(1)]
    	\item $\sum_i x_i$ is a weakly unconditionally Cauchy series (wuC).
    	\item \emph{$S_{S_{\theta}}(\sum_i x_i)$} is a complete space.
        \item \emph{$c_0\subset S_{S_{\theta}}(\sum_i x_i).$}
    \end{enumerate}
\end{theorem}

{\it Proof.} (1)$\Rightarrow$(2): Since $\sum x_i$ is wuC, the following supremum is finite:
$$ H=\sup\left\{\left\|\sum_{i=1}^n a_i x_i\right\|:|a_i|\leq 1,1\leq i\leq n,n\in\mathbb{N}\right\}<+\infty.$$

Let $(a^m)_m\subset S_{S_{\theta}}(\sum_i x_i)$ such that $\displaystyle\lim_m\|a^m-a^0\|_{\infty}=0$, with $a^0\in\ell_{\infty}$. We will prove that $a^0\in S_{S_{\theta}}(\sum_i x_i)$. Let us suppose without any loss of generality that $\|a^0\|_{\infty}\leq 1$. Then, the partial sums $S_k^0=\sum_{i=1}^k a_i^0 x_i$ satisfy $\|S_k^0\|\leq H$ for every $k\in\mathbb{N}$, that is, the sequence $(S_k^0)$ is bounded. Then, $a^0\in S_{S_{\theta}}(\sum_i x_i)$ if and only if $(S_k^0)$ is $S_{\theta}$-summable to some $L\in X$.
According to Theorem \ref{FridyCauchy}, $(S_k^0)$ is lacunary statistically convergent to $L\in X$ if and only if $(S_k^0)$ is a lacunary statistically Cauchy sequence.

Given $\varepsilon>0$ and $n\in\mathbb{N}$,   we obtain statement (2) if we show that there exists a sub-sequence $(S_{k'(r)})$ such that $k'(r)\in I_r$ for every $r$, $\displaystyle\lim_{r\to\infty} S_{k'(r)}=L$ and
$$ d_{\theta}\bigl(\{k\in I_r: ||S_k^0-S^0_{k'(r)}||<\varepsilon\}\bigr)=1.$$

Since $a^m\to a^0$ in $\ell_{\infty}$, there exists $m_0>n$ such that $\|a^m-a^0\|_\infty<\displaystyle\frac{\varepsilon}{4H}$ for all $m>m_0$, and since $S_k^{m_0}$ is $S_{\theta}-$Cauchy, there exists $k'(r)\in I_r $ such that $\displaystyle\lim_{r\to\infty}S_{k'(r)}^{m_0}=L$ for some $L$ and
$$ d_{\theta}\left(\left\{k\in I_r:\|S_{k}^{m_0}-S_{k'(r)}^{m_0}\|<\frac{\varepsilon}{2}\right\}\right)=1.$$

Consider $r\in\mathbb{N}$ and fix $k\in I_r$ such that
\begin{align}\label{eq1}
\displaystyle\|S_k^{m_0}-S_{k'(r)}^{m_0}\|<\frac{\varepsilon}{2}.
\end{align}
We will show that $\|S_k^0-S_{k'(r)}^0\|<\varepsilon$, and this will prove that
$$\left\{k\in I_r:\|S_k^{m_0}-S_{k'(r)}^{m_0}\|<\frac{\varepsilon}{2}\right\}\subset\{k\in I_r:\|S_k^0-S_{k'(r)}^0\|<\varepsilon\}.$$
Since the first set has density 1, the second will also have density 1 and we will be done.

Let us observe first that for every $j\in\mathbb{N}$,
$$\left\|\sum_{i=1}^j\frac{4H}{\varepsilon}(a_i^m-a_i^{m_0})x_i\right\|\leq H,$$
for every $m>m_0$, therefore
\begin{align}\label{eq2}
	\bigl\|S_j^0-S_j^{m_0}\bigr\|=\left\|\sum_{i=1}^j(a_i^0-a_i^{m_0})x_i\right\|\leq\frac{\varepsilon}{4}.
\end{align}
Then, by applying the triangular inequality,
\begin{align*}
	\bigl\|S_k^0-S_{k'(r)}^0\bigr\| & \leq \bigl\|S_k^0-S_k^{m_0}\bigr\|+\bigl\|S_k^{m_0}-S_{k'(r)}^{m_0}\bigr\|+\bigl\|S_{k'(r)}^{m_0}-S_{k'(r)}^{0}\bigr\|\\
    & <\frac{\varepsilon}{4}+\frac{\varepsilon}{2}+\frac{\varepsilon}{4}=\varepsilon.
\end{align*}
where the last inequality follows by applying \eqref{eq1} and \eqref{eq2}, which yields  the desired result.

$(2)\Rightarrow(3)$: Let us observe that if $S_{S_{\theta}}(\sum_i x_i)$ is a complete space, it contains the space of eventually zero sequences $c_{00}$ and therefore the thesis comes, since the supremum norm completion of  $c_{00}$ is $c_0$.

$(3)\Rightarrow (1)$: By way of contradiction, suppose that  the series $\sum x_i$ is not wuC. Therefore there exists $f\in X^{\ast}$ such that $\displaystyle\sum_{i=1}^{\infty}|f(x_i)|=+\infty$. 

CLAIM: We can construct inductively  a sequence $(a_i)_i\in c_0$ such that $$\sum_i a_i f(x_i)=+\infty$$ and $$a_i f(x_i)\geq 0.$$

PROOF: Since $\sum_{i=1}^{\infty}|f(x_i)|=+\infty$, there exists $m_1$ such that $\sum_{i=1}^{m_1}|f(x_i)|>2\cdot2$. 

We define $a_i=\frac{1}{2}$ if $f(x_i)\geq 0$ and $a_i=-\frac{1}{2}$ if $f(x_i)<0$ for $i\in\{1,2,\dots,m_1\}$. 

This implies that $\sum_{i=1}^{m_1}a_i f(x_i)>2$ and $a_i f(x_i)\geq 0$ if $i\in\{1,2,\dots,m_1\}$.

Let $m_2>m_1$ be such that $\sum_{i=m_1+1}^{m_2}|f(x_i)|>2^2\cdot 2^2$. 

We define $a_i=\frac{1}{2^2}$ if $f(x_i)\geq0$ and $a_i=-\frac{1}{2^2}$ if $f(x_i)<0$ for $i\in\{m_1+1,\dots,m_2\}$. Then, $\sum_{i=m_1+1}^{m_2}a_i f(x_i)>2^2$ and $a_i f(x_i)\geq 0$ if $i\in\{m_1+1,\dots,m_2\}$.

So  we have obtained  a sequence $(a_i)_i\in c_0$ with the above properties.
\vskip5mm
Now  we will prove that the sequence $S_k=\sum_{i=1}^k a_i f(x_i)$ is not $S_{\theta}$-summable to any $L\in\R$. 
By way of contradiction, suppose that it is  $S_{\theta}$-summable to $L\in \R$, then we have
$$\frac{1}{h_r}\card(\{k\in I_r:|S_k-L|\geq\varepsilon\})=\frac{1}{h_r}\sum_{\substack{k=k_{r-1}\\ |S_k-L|\geq\varepsilon}}^{k_r}1\underset{r\to\infty}{\to} 0.
$$
Since $S_k$ is an increasing sequence and $S_k\to\infty$, there exists $k_0$ such that $|S_k-L|\geq\varepsilon$ for every $k\geq k_0$. Let us suppose that $k_r>k_0$ for every $r$. Hence,
$$\frac{1}{h_r}\sum_{\substack{k=k_{r-1}\\ |S_k-L|\geq\varepsilon}}^{k_r}1=\frac{h_r}{h_r}=1\underset{r\to\infty}{\nrightarrow} 0,$$
which is a contradiction. This implies that $S_k$ is not $S_\theta$-convergent and this is a contradiction with (3).

\section{The strong lacunary summability space}

Let $\sum_i x_i$ be a series in a real Banach space $X$ and $\theta=(k_r)$ a lacunary sequence. We define
$$S_{N_{\theta}}\left(\sum_i x_i\right)=\left\{(a_i)_i\in\ell_{\infty}: \sum_{i} a_i x_i\mbox{ is }N_{\theta}\mbox{-summable}\right\}$$
endowed with the supremum norm. This space will be called the space of $N_{\theta}$-summability associated to the series $\sum_i x_i$. The following theorem characterizes the completeness of the space $S_{N_{\theta}}\bigl(\sum_i x_i\bigr)$.

\begin{theorem}\label{N_tcompleto}
	Let $X$ be a real Banach space and $\theta=(k_r)$ a lacunary sequence. The following conditions are equivalent:
    \medskip
    \begin{enumerate}[(1)]
    	\item $\sum_i x_i$ is a weakly unconditionally Cauchy series (wuC).
    	\item \emph{$S_{N_{\theta}}(\sum_i x_i)$} is a complete space.
        \item \emph{$c_0\subset S_{N_{\theta}}(\sum_i x_i).$}
    \end{enumerate}
\end{theorem}

{\it Proof.} (1)$\Rightarrow$(2): Since $\sum x_i$ is wuC, the following supremum is finite
$$H=\sup\left\{\left\|\sum_{i=1}^n a_i x_i\right\|:|a_i|\leq 1,1\leq i\leq n,n\in\mathbb{N}\right\}<+\infty.$$

Let $(a^m)_m\subset S_{N_{\theta}}(\sum_i x_i)$ such that $\displaystyle\lim_m\|a^m-a^0\|_{\infty}=0$, with $a^0\in\ell_{\infty}$. 

We will prove that $a^0\in S_{N_{\theta}}(\sum_i x_i)$. 

Without  loss of generality  we can  suppose that $\|a^0\|_{\infty}\leq 1$. Therefore the partial sums $S_k^0=\sum_{i=1}^k a_i^0 x_i$ satisfy $\|S_k^0\|\leq H$ for every $k\in\mathbb{N}$, that is, the sequence $(S_k^0)$ is bounded. Hence  $a^0\in S_{N_{\theta}}(\sum_i x_i)$ if and only if $(S_k^0)$ is 
$N_{\theta}$-summable to some $L\in X$.  Since $(S_k^0)$ is bounded, it is sufficient to show that $(S_k)$ is $S_{\theta}$-convergent, thanks to to Fridy and Orhan's Theorem \cite[Theorem 2.1]{Fridy02} (see Theorem \ref{FridyAcotadas}). The result follows  analogously  as in Theorem \ref{S_tcompleto}.

$(2)\Rightarrow(3)$: It is sufficient to notice that  $S_{S_{\theta}}(\sum_i x_i)$ is a complete space and it contains the space of eventually zero sequences $c_{00}$, so it  contains the  completion of $c_{00}$ with respect to the supremum norm, hence it contains $c_0$.

$(3)\Rightarrow (1)$: By way of contradiction, suppose that  the series $\sum x_i$ is not wuC.  Therefore there exists $f\in X^{\ast}$ such that $\displaystyle\sum_{i=1}^{\infty}|f(x_i)|=+\infty$. We can  construct inductively  a sequence $(a_i)_i\in c_0$ as in Theorem \ref{S_tcompleto} such that $\sum_i a_i f(x_i)=+\infty$ and $a_i f(x_i)\geq 0$.

The sequence $S_k=\sum_{i=1}^k a_i f(x_i)$ is not $N_{\theta}$-summable to any $L\in\R$. 

As $S_k\to \infty$, for every  $A>0$, there exists $k_0$ such that $|S_k|>A$ if $k\geq k_0$. Then we have
$$\frac{1}{h_r}\sum_{k\in I_r}|S_k|>\frac{h_r A}{h_r}=A.$$
Hence $S_k$ is not $N_{\theta}$-summable to any $L\in \R$, otherwise 
$$\infty\leftarrow\frac{1}{h_r}\sum_{k\in I_r}|S_k|\leq |L| + \frac{1}{h_r}\sum_{k\in I_r}|S_k-L|\to |L|$$

We can conclude that $S_k$ is not $N_\theta$-convergent, a contradiction with (3).$\qed$

\section{Characterizations of the completeness of a Banach space}

A Banach space $X$ can be characterized by  the completeness of the space $S_{N_{\theta}}(\sum_i x_i)$ for every wuC series $\sum_i x_i$, as we
will show next.

\begin{theorem}\label{teor3.5}
Let $X$ be a normed real vector space. Then $X$ is complete if and only if
\emph{$S_{N_{\theta}}(\sum_i x_i)$} is a complete space for every weakly unconditionally Cauchy series (wuC) $\sum_i x_i$.
\end{theorem}
{\it Proof.} Thanks to  Theorem \ref{S_tcompleto}, the condition is necessary. 

Now suppose that  $X$ is not complete, hence there exists 
a series $\sum x_i$ in $X$ such that $\displaystyle\|x_i\|\leq \frac{1}{i2^i }$ and $\sum x_i=x^{\ast\ast}\in X^{\ast\ast}\setminus X$. 

We will construct a wuC series $\sum_i y_i$ such that $S_{N_{\theta}}(\sum_i y_i)$ is not complete, a contradiction. 

Set $\displaystyle S_N=\sum_{i=1}^N x_i$. As $X^{\ast\ast}$ is a Banach space endowed with the dual topology, $\disp\sup_{\|y^\ast\|\leq 1}|y^\ast(S_N)-x^{\ast\ast}(y^\ast)|$ tends to $0$ as $N\to\infty$, that is,
\begin{align}\label{(A)}
\lim_{N\to+\infty}y^\ast(S_N)=\lim_{N\to+\infty}\sum_{i=1}^N y^\ast(x_i)=x^{\ast\ast}(y^\ast), \mbox{ for every } \|y^\ast\|\leq 1.
\end{align}
Put $y_i=i x_i$ and let us observe that $\|y_i\|<\frac{1}{2^i}$. Therefore  $\sum y_i$ is absolutely convergent, thus it is unconditionally convergent and weakly unconditionally Cauchy. 

We claim that the series $\displaystyle\sum_i \frac{1}{i}y_i$ is not $N_{\theta}$-summable in $X$.

By way of contradiction suppose that $S_N=\sum_{i=1}^N \frac{1}{i}y_i$ is $N_{\theta}$-summable in $X$, i.e., there exists $L\in X$ such that $\disp\lim_{r\to\infty}\frac{1}{h_r}\sum_{i\in I_r}\|S_i-L\|=0$. This implies that
\begin{align}\label{(B)}
\lim_{r\to+\infty}\frac{1}{h_r}\sum_{i\in I_r} y^\ast(S_i)=y^\ast(L), \mbox{ for every } \|y^\ast\|\leq 1.
\end{align}

From equations \eqref{(A)} and \eqref{(B)}, the uniqueness of the limit and since $N_\theta$ is a regular method, we have  $x^{\ast\ast}(y^\ast)=y^\ast(L)$ for every $\|y^\ast\|\leq 1$, so we obtain $x^{\ast\ast}=L\in X$, a contradiction. Hence $S_N=\sum_{i=1}^N \frac{1}{i}y_i$ is not $N_{\theta}$-summable to any $L\in X$.

Finally, let us observe that, since $\sum_i y_i$ is a weakly unconditionally Cauchy series and  $S_N=\sum_{i=1}^N\frac{1}{i}y_i$ is not $N_{\theta}$-summable, we have  $(\frac{1}{i})\notin {S_{N_{\theta}}(\sum_i y_i)}$ and this means that $c_0\nsubseteq S_{N_{\theta}}(\sum_i y_i)$ which is a contradiction with  Theorem \ref{N_tcompleto}(3), so   the proof is complete.$\qed$
\vskip5mm

By  a similar argument and taking into account Theorem \ref{FridyAcotadas}, we  have also the characterization for the $S_\theta$-summability:

\begin{theorem}\label{Completeness_St}
Let $X$ be a normed real vector space. Then $X$ is complete if and only if
\emph{$S_{S_{\theta}}(\sum_i x_i)$} is a complete space for every weakly unconditionally Cauchy series (wuC) $\sum_i x_i$.
\end{theorem}

 Let $0<p<+\infty$, the sequence  $(x_n)$ is said to be strongly p-Ces\`aro or $w_p$-summable if there is $L\in X$ such that $$\lim_n \frac{1}{n}\sum_{i=1}^{n} \Vert x_i-L\Vert^p=0;$$ in this case
we will write $(x_k)\rightarrow_{w_p} L$ and $L = w_p-\lim_n x_n.$ 
Let $\sum x_i$ be a series in a real Banach space $X$, let us define $$
	S_{w_p}\left(\sum_i x_i\right)=\left\{(a_i)_i\in\ell_{\infty}: \sum_{i} a_i x_i\;\mbox{ is} \;w_p\mbox{-summable}\right\}$$
endowed with the supremum norm.

We refer to \cite{SalaFilomat} for other properties  of  the space \emph{$S_{w_p}(\sum_i x_i)$}.

Finally, from Theorem \ref{teor3.5}, Theorem \ref{Completeness_St} and \cite[Theorem 3.5]{SalaFilomat}, we derive the following corollary. 

\begin{corollary}
Let $X$ be a normed real vector space and $p\geq 1$. Then the following items are equivalent:
\begin{enumerate}
\item X is complete.
\item $S_{N_{\theta}}(\sum_i x_i)$ is a complete space for every weakly unconditionally Cauchy series (wuC) $\sum_i x_i$.
\item $S_{S_{\theta}}(\sum_i x_i)$ is a complete space for every weakly unconditionally Cauchy series (wuC) $\sum_i x_i$.
\item $S_{\emph{w}_p}(\sum_i x_i)$ is a complete space for every weakly unconditionally Cauchy series (wuC) $\sum_i x_i$.
\end{enumerate}
\end{corollary}

\section*{Acknowledgements}

This work is supported by the FQM-257 research group of the University of Cádiz and the Research Grant PGC-101514-B-100 awarded by the Spanish Ministry of Science, Innovation and Universities and partially funded by the European Regional Development Fund.


\begin{thebibliography}{10}

\bibitem{boos}
J.~Boos.
\newblock {\em Classical and modern methods in summability}.
\newblock Oxford Mathematical Monographs. Oxford University Press, Oxford,
  2000.
\newblock Assisted by Peter Cass, Oxford Science Publications.

\bibitem{connorganichevkadets}
J.~Connor, M.~Ganichev, and V.~Kadets.
\newblock A characterization of {B}anach spaces with separable duals via weak
  statistical convergence.
\newblock {\em J. Math. Anal. Appl.}, 244(1):251--261, 2000.

\bibitem{Diestel1984}
J.~Diestel.
\newblock {\em Sequences and Series in Banach Spaces}.
\newblock Graduate Texts in Mathematics. Springer New York, 1984.

\bibitem{fekete}
M.~Fekete.
\newblock Viszg{\'a}latok a fourier-sorokr{\'o}l (research on fourier series).
\newblock {\em Math. {\'e}s term{\'e}sz}, 34:759--786, 1916.

\bibitem{Freedman78}
A.~R. Freedman, J.~J. Sember, and M.~Raphael.
\newblock Some {C}es\`aro-type summability spaces.
\newblock {\em Proc. London Math. Soc. (3)}, 37(3):508--520, 1978.

\bibitem{Fridy02}
J.~A. Fridy and C.~Orhan.
\newblock Lacunary statistical convergence.
\newblock {\em Pacific J. Math.}, 160(1):43--51, 1993.

\bibitem{Fridy01}
J.~A. Fridy and C.~Orhan.
\newblock Lacunary statistical summability.
\newblock {\em J. Math. Anal. Appl.}, 173(2):497--504, 1993.

\bibitem{hardy}
G.~H. Hardy and J.~E. Littlewood.
\newblock Sur la s{\'e}rie de fourier d'une fonction {\'a} carr{\'e} sommable.
\newblock {\em CR Acad. Sci. Paris}, 156:1307--1309, 1913.

\bibitem{JIneq3}
A.~Kilicman and S.~Borgohain.
\newblock Some new lacunary statistical convergence with ideals.
\newblock {\em J. Inequal. Appl.}, pages Paper No. 15, 9, 2017.

\bibitem{kolk}
E.~Kolk.
\newblock The statistical convergence in {B}anach spaces.
\newblock {\em Tartu \"Ul. Toimetised}, (928):41--52, 1991.

\bibitem{SalaFilomat}
F.~Le\'on-Saavedra, S.~Moreno-Pulido, and A.~Sala-P\'erez.
\newblock Completeness of a normed space via strong $p$-{C}es\`aro summability.
\newblock {\em Filomat}, 33(10):3013--3022, 2019.

\bibitem{Leon-Saavedra2019}
F.~León-Saavedra, F.~J. Pérez-Fernández, M.~P. Romero de~la Rosa, and
  A.~Sala.
\newblock Ideal convergence and completeness of a normed space.
\newblock {\em Mathematics}, 7(10):1--11, 2019.

\bibitem{JIneq5}
S.~A. Mohiuddine, A.~Alotaibi, and M.~Mursaleen.
\newblock A new variant of statistical convergence.
\newblock {\em J. Inequal. Appl.}, 2013.

\bibitem{MursaleenLibro}
M.~Mursaleen.
\newblock {\em Applied summability methods}.
\newblock SpringerBriefs in Mathematics. Springer, Cham, 2014.

\bibitem{Aizpuru2000}
F.~J. P\'erez-Fern\'andez, F.~Ben\'{i}tez-Trujillo, and A.~Aizpuru.
\newblock Characterizations of completeness of normed spaces through weakly
  unconditionally {C}auchy series.
\newblock {\em Czechoslovak Math. J.}, 50(125)(4):889--896, 2000.

\bibitem{zeller-beekmann}
K.~Zeller and W.~Beekmann.
\newblock {\em Theorie der {L}imitierungsverfahren}.
\newblock Zweite, erweiterte und verbesserte Auflage. Ergebnisse der Mathematik
  und ihrer Grenzgebiete, Band 15. Springer-Verlag, Berlin-New York, 1970.

\bibitem{zygmund}
A.~Zygmund.
\newblock {\em Trigonometrical series}.
\newblock Dover Publications, New York, 1955.

\end{thebibliography}
\end{document}